\documentclass[preprint, 12pt]{elsarticle}
\usepackage{amsthm}
\usepackage{amssymb,amsfonts}
\usepackage[tbtags,fleqn]{amsmath} 
\allowdisplaybreaks  

\usepackage{caption}
 \newtheorem{thm}{Theorem}

\newdefinition{rmk}{Remark}
\newproof{pf}{Proof}
\newproof{pot}{Proof of Theorem \ref{thm2}}

\makeatletter
\renewcommand\@biblabel[1]{}
\renewenvironment{thebibliography}[1]
     {\section*{\refname}%
      \@mkboth{\MakeUppercase\refname}{\MakeUppercase\refname}%
      \list{\@biblabel{\@arabic\c@enumiv}}%
           {\settowidth\labelwidth{\@biblabel{#1}}%
            \leftmargin\labelwidth
            \advance\leftmargin\labelsep
            \advance\leftmargin by 2em%
            \itemindent -2em%
            \@openbib@code
            \usecounter{enumiv}%
            \let\p@enumiv\@empty
            \renewcommand\theenumiv{\@arabic\c@enumiv}}%
      \sloppy
      \clubpenalty4000
      \@clubpenalty \clubpenalty
      \widowpenalty4000%
      \sfcode`\.\@m}
     {\def\@noitemerr
       {\@latex@warning{Empty `thebibliography' environment}}%
      \endlist}
\makeatother

\begin{document}

\begin{frontmatter}

\title{Analysis of the Optimal Resource Allocation for a Tandem Queueing System}

\author{Liu Zaiming}
\ead{math\_lzm@csu.edu.cn}
\author{Chen Gang}
\ead{chengmathcsu@163.com}
\author{Wu Jinbiao$^*$ }
\ead{$^*$ Corresponding author: wujinbiao@csu.edu.cn}

\address{School of Mathematics and Statistics, Central
South University, Changsha 410083, Hunan, PR China}

\begin{abstract}
 In this paper, we study a controllable tandem queueing system consisting of two nodes and a controller, in which customers arrive according to a Poisson process and must receive service at both nodes before leaving the system. A decision maker dynamically allocates the number of service resource
 to each node facility according to the number of customers in each node. In the model, the objective is to minimize the long-run average costs. We cast these problems as Markov decision problems by dynamic programming approach and derive the monotonicity of the optimal allocation policy and the relationship between the two nodes' optimal policy. Furthermore, we get the conditions under which the optimal policy is unique and has the bang-bang control policy property.
\end{abstract}

\begin{keyword}

Markov decision problem \sep Tandem system\sep Optimal policy \sep Dynamic programming
\sep Average costs

\end{keyword}

\end{frontmatter}

\section{Introduction}\label{sec1}

We consider a controllable tandem queueing system consisting of two nodes and a controller. A decision maker can assign a number of service resource to each node. The study of the controllable tandem queueing system
 is motivated by its wide applications in manufacturing, computer systems, voice and data communications, and vehicular traffic flow. The theory of controllable queueing systems has often been studied for optimal control of admission, servicing, dynamic pricing, routing and scheduling of jobs in queues or networks of queues. These works are discussed in Stidham and Weber (1993), Yang et al. (2011) and {\c{C}}il et al. (2011). The controllable queueing systems based on the theory of Markov, semi-Markov and regenerative decision processes can be found in Morozov and Steyaert (2013). Using the theory of the queueing system, we often cast the optimal problems as Markov decision problems (MDP). In order to get the properties of the optimal policy, the properties (such as the monotonicity, convexity property) of relative value function (when we consider the long-run average criteria) should be first considered. The key of the method is dynamic programming. For more details, we can see the paper written by Koole (1998) and {\c{C}}il et al. (2009).

Based on the application background, the problems of the service resource control in different queueing systems have been investigated. Rykov and Efrosinin (2004) considered a multi-server controllable queueing system with heterogeneous servers, and several monotonicity properties of optimal policies are proved. Iravani et al. (2007) studied the optimal service scheduling in nonpreemptive finite-population queueing systems. The single-queue systems of the optimal resource allocation policy were considered by Yang et al. (2013). Efrosinin et al. (2014) analyzed a tandem queueing system of admission optimal policy.

Of particular relation to the present work are the works of Rosberg et al. (1982) and Ahn et al. (2002) where only the customer's holding cost was considered. Rosberg et al. (1982) considered the optimal control of service in tandem queues where the service rate in node 1 can be selected from a compact set and constant in node 2. Optimal control of a two-stage tandem queues system with flexible servers was discussed in Ahn et al. (2002) where only two flexible servers were considered under two different scenarios and they obtained the exhaustive optimal policy. Kaufman et al. (2005) considered the problem on the agile, temporary workforce into a tandem queueing system in which the relationship between the service rate and the number of the service resource is linear and the service resource costs in different nodes have the same cost function. However, different from the previous studies about resource allocation control problem, the two nodes in our model have the different holding cost rate and service resource cost function in the objective (long-run average cost). The main contribution of this paper is that we derive the monotonicity of the optimal allocation policy and the relationship between the two nodes' optimal policy. Furthermore, we get the conditions under which the optimal policy is unique and the bang-bang control policy is established.

The rest of the paper is organized as follows. In Section \ref{sec2}, the model is formulated in detail based on the controllable Markov decision problem. The characteristics of the optimization problem and the optimality equation are derived in Section \ref{sec3}.
In Section \ref{sec4}, structural properties of the optimal policy and main results of the paper are given. Finally, some further discussions and conclusions are given in Section \ref{sec5}.

 \section{Model Description }\label{sec2}
We consider a tandem queueing system with two nodes. Customers arrive at node 1 from outside the system according to a Poisson process with parameter $\lambda$ and have exponentially distributed service requirement at each node. After receiving service at node 1, customers proceed immediately to node 2 and receive service before leaving the system. A decision maker can assign a number of service resource to each node. The service rate of a customer depends on the number of service resource assigned to the customer precisely. When a customer has been allocated $a$ server resources, the service duration of that customer in node $i$ is exponentially distributed with parameter $\mu_{i}(a),i=1,2$, which is strictly increasing in $a$. Without loss of generality, we assume that $\mu_{i}(0)=0,i=1,2$. At any decision epoch, the decision maker decides to choose the number of server resources to node 1 from a compact set $A=[0,a_{max}]$, and to node 2 from a compact set $B=[0,b_{max}]$ at the same time. Each node has a single infinite-size FCFS queue. The interarrival and service times are assumed to be mutually independent. We assume that the stability condition $\lambda<\mu_{1}(a_{max}),\lambda<\mu_{2}(b_{max})$ holds. Figure 1 gives an illustration of the system.

\begin{figure}[htp]
\centering
\includegraphics[width=4.3in]{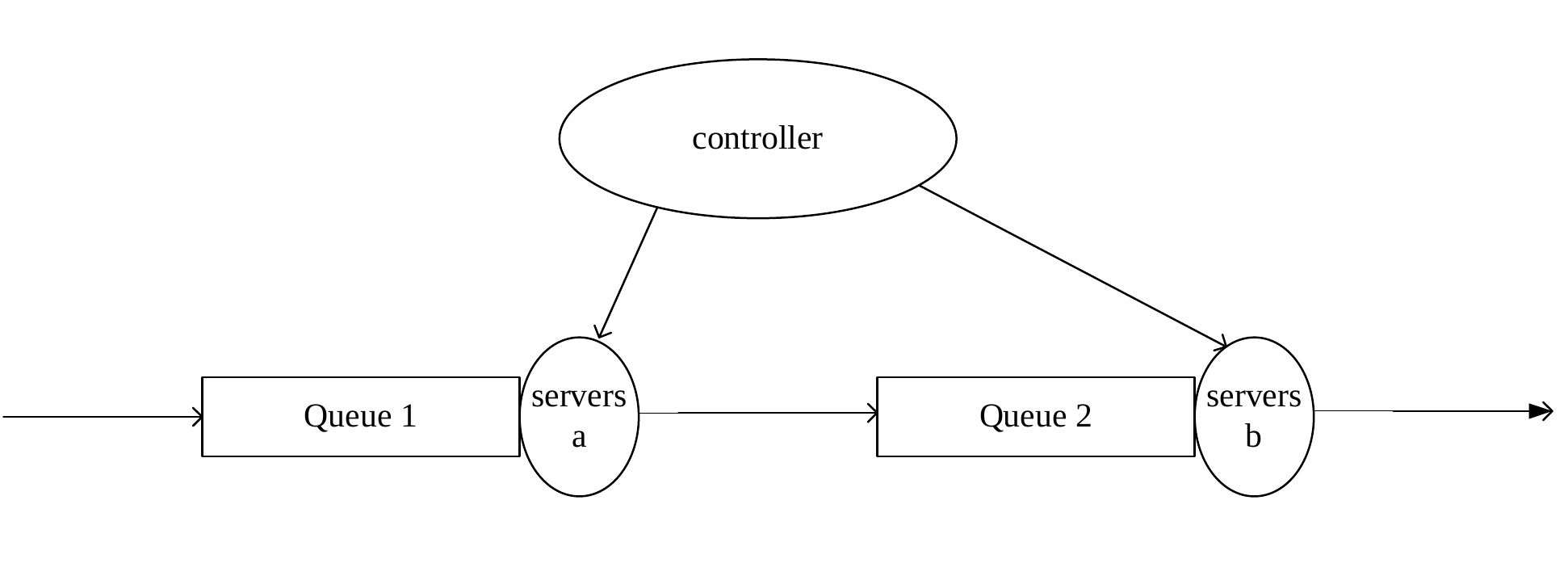} 
\caption {The controllable tandem queueing systems}
\label{Fig.1}
\end{figure}

We consider the following cost structure in the system. Our objective is to obtain dynamic management policy that minimizes the long-run average costs.

(1) resources cost: when the node $i$ uses $a$ resources, a cost of $c_{i}(a),i=1,2$ is incurred by the system per unit time (here $c_{i}(a)$ is a continuous function and strictly increasing in $a$. Without loss of generality, we assume that $c_{i}(0)=0,i=1,2$).

(2) holding cost: holding costs are incurred at rate $h_{1}$ and $h_{2}$ per unit time for each customer in node 1 and 2, respectively.

Let $X_{i}(t)$ denote the number of customers at node $i,i=1,2$. The system evolves as a continuous-time Markov process
\begin{equation*}\label{eq:nc1}
 \{X(t),t\geq0\}=\{(X_{1}(t),X_{2}(t)),t\geq0\}.
\end{equation*}
The notations $l_{i}(x),i=1,2$, will be used to specify the certain components of the vector state $x\in E$.

The system state space is: $E={x=(x_{1},x_{2})\in N^{2}}$, with $N={0,1,2,...}$.

It is assumed that the model is stable and conservative. The transition rate under a control action $(a,b)$ is given by\\
\[Q_{xy}(a,b)=\left\{
\begin{array}{cc}
\lambda & y=x+e_{1};\\
\mu_{1}(a) &\    \ \ y=x-e_{1}+e_{2},l_{1}(x)>0;\\
\mu_{2}(b) & y=x-e_{2},l_{2}(x)>0;\\
0  & \mbox{else},
\end{array}
\right.
\]
where
\begin{equation*}\label{eq:nc1}
 Q_{xy}(a,b)\geq0, y\neq x, Q_{xx}(a,b)=-Q_{x}(a,b)=-\sum_{y\neq x} Q_{xy}(a,b), Q_{x}(a,b)<\infty.
\end{equation*}
Here $e_{i}$ is the 2-dimensional vector with 1 in the  $i$th coordinate and 0 elsewhere, $i=1,2$.

The problem of the decision maker is to derive an optimal policy based on the number of customers in each node that minimizes the long-run average costs. We cast the customer resource management problem as a Markov decision problem. The set of decision epochs corresponds to the set of all arrivals, service completions, and dummy transitions due to uniformization. The controllable system associated with a Markov process is a five-tuple

$\{E,D=(A,B),Q(f),c_{i}(a),h_{i}\},i=1,2$,\\
in which $Q(f)$ is the transition matrix of the queueing system under the policy $f$.

We consider the stationary Markov policy $f:E\rightarrow D$ with $f=(f_{1},f_{2})$. Due to the Markov property, it is clear that the optimal policy depends only on the current state regardless of $t$. More precisely, when the system state is $x=(x_{1},x_{2})$, the controller makes an action $f_{1}(x_{1})=a\in A,f_{2}(x_{2})=b\in B$. The action of the service resource to node $i$ only depends on the current number of customers in node $i$.

\section{Optimization problem and optimality equation}\label{sec3}
For every fixed stationary policy $f$, we assume that the process $\{X(t),t\geq0\}$ with state space $E$ is an irreducible, positive recurrent Markov process. As it is known from Tijms (1994), for ergodic Markov process with the long-run average cost per unit of time for the policy $f$ coincides with corresponding assemble average,
\begin{equation}\label{eq:nc1}
g(f)=\lim_{t\rightarrow\infty}u(x,t)^{f}/t=\sum_{i=1}\sum_{j=1}[c_{1}(f_{1}(i))+c_{2}(f_{2}(j))+h_{1}i+h_{2}j]\pi_{ij}(f),
\end{equation}
in which $u(x,t)^{f}$ denotes the total expected costs up to time $t$ when the system starts in state $x$ and $\pi_{ij}(f)$ denotes a stationary probability of the process  under policy $f$. The goal is to find a policy $f^{*}$ that minimizes the long-term average costs:
\begin{equation}\label{eq:nc1}
g(f^{*})=\min_{f}g(f).
\end{equation}
In order to find the optimal policy $f^{*}$ that minimizes the total average cost, we construct a discrete-time equivalent of the original system by using the standard tools of uniformization and normalization. Without loss of generality, we assume that $\lambda+\mu_{1}(a_{max})+\mu_{2}(b_{max})=1$. Now we consider a real-valued function $v(x)$ that plays the role of the relative value function, i.e., the asymptotic difference in total costs that results from starting the process in state $x$ instead of some reference state. As it is well known, the optimal policy $f$ and the optimal average cost $g$ are the solutions of the optimality equation
 \begin{equation*}\label{eq:nc1}
Tv(x)=v(x)+g,
\end{equation*}
where $T$ is the dynamic programming operator acting on $v$, defined as follows
\begin{equation}\label{eq:nc1}
Tv(x)=\lambda v(x+e_{1})+\Sigma_{i=1,2}T_{i}v(x)+\Sigma_{i=1,2}h_{i}l_{i}(x),
\end{equation}
here
\begin{eqnarray}\label{eq:nc1}
T_{1}v(x)=\min_{a\in A}\{\mu_{1}(a)v(x-e_{1}+e_{2})+[\mu_{1}(a_{max})-\mu_{1}(a)]v(x)+c_{1}(a)\},
\end{eqnarray}
\begin{eqnarray}\label{eq:nc1}
T_{2}v(x)=\min_{b\in B}\{\mu_{2}(b)v(x-e_{2})+[\mu_{2}(b_{max})-\mu_{2}(b)]v(x)+c_{2}(b)\}.
\end{eqnarray}

The first term in the expression $Tv(x)$ models the arrivals of customers to node 1 from outside the system and the last one the customer holding cost. Similarly the first term in the expression $T_{1}v(x)$ corresponds to a customer who finished his service in node 1 and into node 2 and the second one the uniformization constant. The last one in $T_{1}v(x)$ is the resources cost in node 1. The first term in the expression $T_{2}v(x)$ corresponds to a customer who finished his service in node 2 and the second one the uniformization constant. The last one in $T_{2}v(x)$ is the resources cost in node 2.

According to (1), we can solve another optimization problem: if $c_{i}\equiv0,h_{i}=1,i=1,2$, then (2) is equivalent to minimization of the mean number of customers in the queueing system.

 \section{Structural properties of the optimal policy}\label{sec4}
In this section, we focus on deriving the optimal policy. However, the optimal policy possesses structural properties that provide fundamental insight, and this also enables one to determine the optimal policy with less computational effort due to a reduction of the solution search space.

In order to study the structure, in principle, one needs to solve the optimal equation $Tv(x)=v(x)+g$. However it is hard to solve analytically in practice. It can be obtained by recursively defining $v_{n+1}=Tv_{n}$ for arbitrary $v_{0}$. We know that the actions converge to the optimal policy as $n\rightarrow\infty$. For existence and convergence of the solutions and optimal policy we refer to Aviv and Federgruen (1999) and Sennott (2009).
The backward recursion equation is given by
 \begin{equation*}\label{eq:nc1}
v_{n+1}(x)=\lambda v_{n}(x+e_{1})+\sum_{i=1,2}T_{i}v_{n}(x)+\sum_{i=1,2}h_{i}l_{i}(x).
\end{equation*}

For ease of notation, we define the set of the optimal policy in state $x$ by:
 \begin{equation*}\label{eq:nc1}
f(x)=(f_{1}(x_{1}),f_{2}(x_{2})) \\f_{1}(x_{1})=argT_{1}v(x) \\ f_{2}(x_{2})=argT_{2}v(x).
\end{equation*}

By using the optimality equation, we can get the properties of relative value function as follows:
\section*{Property 4.1 (non-decreasingness)}
(i) $ v(x+e_{i})\geq v(x),i=1,2$ for all $x\in E$,\\

(ii) if $2h_{2}\geq h_{1}$ then $v(x-e_{1}+e_{2})\geq v(x-e_{2})$ for all $x=(x_{1},x_{2})\in E$ and $x_{1}\geq1,x_{2}\geq1$,\\

(iii) if $h_{1}\geq h_{2}$ then $v(x)\geq v(x-e_{1}+e_{2})$ for all $x=(x_{1},x_{2})\in E$ and $x_{1}\geq1,x_{2}\geq1$.

\section*{Property 4.2 (quasi-convexity)}

(i) $ v(x+e_{2})-2v(x)+v(x-e_{2})\geq0$, for all $ x=(x_{1},x_{2})\in E$ and $x_{2}\geq1$,

(ii) $ v(x+e_{1}-e_{2})-2v(x)+v(x-e_{1}+e_{2})\geq0$, for all $ x=(x_{1},x_{2})\in E$ and $x_{1}\geq1,x_{2}\geq1$.

Next we show some structure properties of the optimal policy, based on the structure properties of the relative value function above.
\begin{thm}\label{thm1}
The optimal policy has the monotonicity property, i.e.,

(i) if $ b_{1}\in argT_{2}v(x+e_{2}),b_{2}\in argT_{2}v(x)$, then $b_{1}\geq b_{2}$ for all $x=(x_{1},x_{2})\in E$.

(ii) if $ a_{1}\in argT_{1}v(x+e_{1}),a_{2}\in argT_{1}v(x)$, then $a_{1}\geq a_{2}$ for all $x=(x_{1},x_{2})\in E$.\\
\end{thm}
The proof of Property 4.1 is given in Appendix A. The proof of Property 4.2 and Theorem 1 are given in Appendix B.\\

Based on Property 4.1, we give the relationship between the two nodes' optimal policy under some conditions.

\begin{thm}\label{thm3}
Assume that $c_{1}(a)-c_{1}(b)\geq c_{2}(a)-c_{2}(b)$ and $\mu_{2}(a)-\mu_{2}(b)\geq \mu_{1}(a)-\mu_{1}(b)$ when $a\geq b$. Then if $a\in argT_{1}v(x),b\in argT_{2}v(x)$, we have $b\geq a$ for all $x=(x_{1},x_{2})\in E$ and $x_{1}\geq1,x_{2}\geq1$.
\end{thm}

\begin{pf}
Let $(a\in argT_{1}v(x),b\in argT_{2}v(x))$ be an arbitrary optimal policy for node 1 and 2 in state $x$, respectively. The proof is by contradiction. Suppose that $b<a$, then we compare the policy $(a,b)$ with the policy $(b,a)$. We have:
\begin{eqnarray*}
\lefteqn{T_{a,b}v_{n}(x)-T_{b,a}v_{n}(x)}\\
&&=[\mu_{1}(a)v(x-e_{1}+e_{2})+[\mu_{1}(a_{max})-\mu_{1}(a)] v(x)+c_{1}(a)]\\
  &&\   \  +[\mu_{2}(b)v(x-e_{2})+[\mu_{2}(b_{max})-\mu_{2}(b)]v(x)+c_{2}(b)]\\
   &&\   \  -[\mu_{1}(b)v(x-e_{1}+e_{2})+[\mu_{1}(b_{max})-\mu_{1}(b)]v(x)+c_{1}(b)]\\
   &&\   \  -[\mu_{2}(a)v(x-e_{2})+[\mu_{2}(a_{max})-\mu_{2}(a)]v(x)+c_{2}(a)]\\
&&=[\mu_{1}(a)-\mu_{1}(b)][v(x-e_{1}+e_{2})-v(x)]-[\mu_{2}(a)-\mu_{2}(b)][v(x-e_{2})-v(x)]\\
   &&\   \  +c_{1}(a)-c_{1}(b)-c_{2}(a)+c_{2}(b)\\
&&\geq [\mu_{1}(a)-\mu_{1}(b)][v(x-e_{1}+e_{2})-v(x-e_{2})]+c_{1}(a)-c_{1}(b)-c_{2}(a)+c_{2}(b)\\
&&\geq 0.
\end{eqnarray*}
The first equality is based on the definition of the operators $T_{1}$ and $T_{2}$. The second equality follows by rearranging the terms. The first inequality follows the condition $\mu_{2}(a)-\mu_{2}(b)\geq \mu_{1}(a)-\mu_{1}(b)$ when $a\geq b$. This implies that $a$ and $b$ is not an optimal policy for node 1 and 2 in state $x$, respectively. Hence, $b\geq a$.

From the above theorem we can conclude that under some conditions the optimal size of the service resources allocate to node 1 is less than that to node 2. We find that the optimal size of the resource allocate to each node depends on the resource cost variation $c(a)-c(b)$ and the service rate variation $\mu(a)-\mu(b)$ in each node.
\end{pf}

We are now ready to give some conditions under which the optimal policy is unique and is a bang-bang control policy.

\begin{thm}\label{thm3}
The following properties hold

(i) if the functions $m _{1}(a)=\frac{c'_{1}(a)}{\mu'_{1}(a)}$ and $m _{2}(b)=\frac{c'_{2}(b)}{\mu'_{2}(b)}$ are monotonous on $a\in A,b\in B$, then the optimal policy is unique.

(ii) $argT_{1}v(0)=\{0\},argT_{2}v(0)=\{0\}$.

(iii) if the functions $\frac{c_{1}(a)}{\mu_{1}(a)}$ and $\frac{c_{2}(b)}{\mu_{2}(b)}$ are non-increasing, $\frac{c'_{1}(a)}{\mu'_{1}(a)}>\frac{c_{1}(a)}{\mu_{1}(a)}$ and $\frac{c'_{2}(b)}{\mu'_{2}(b)}>\frac{c_{2}(b)}{\mu_{2}(b)}$ for all $a\in (0,a_{max}),b\in (0,b_{max})$, then the optimal policy is a bang-bang control policy. i.e., $argT_{1}v(x)=\{0,a_{max}\},argT_{2}v(x)=\{0,b_{max}\}$ for all $x\in E$.

\end{thm}
\begin{pf}
To prove part (i), we consider the optimal policy $a$ in node 1 service resource allocation. In our event operator $T_{1}$ for node 1 defined in equation (3), we have the following minmization problem:
\begin{eqnarray*}\label{eq:nc1}
T_{1}v(x)=\min_{a\in A}\{\mu_{1}(a)v(x-e_{1}+e_{2})+[\mu_{1}(a_{max})-\mu_{1}(a)]v(x)+c_{1}(a)\}.
\end{eqnarray*}
Rearranging the first-order optimality condition of the above problem, we have:
\begin{eqnarray*}\label{eq:nc1}
\frac{c'_{1}(a)}{\mu'_{1}(a)}=v(x)-v(x-e_{1}+e_{2}).
\end{eqnarray*}
Because the allocation resource action $a\in A=[0,a_{max}]$, the optimal policy $a$ must be the solution of the above equation. Since the function $m _{1}(a)=\frac{c'_{1}(a)}{\mu'_{1}(a)}$ is monotonous on $a\in A$, there is a unique $a$ solving the above equation. Hence the optimal policy for node 1 is unique. The part (i) for node 2 can be proved in a similar manner.\\

To prove part (ii), we consider the optimal policy $a$ in node 1 service resource allocation. As the problem is defined in equation (3), we have
\begin{eqnarray*}\label{eq:nc1}
T_{1}v(0)=\min_{a\in A}\{\mu_{1}(a)v(0)+[\mu_{1}(a_{max})-\mu_{1}(a)]v(0)+c_{1}(a)\},
\end{eqnarray*}
which immediately implies that $argT_{1}v(0)=\{0\}$. The part (ii) for node 2 that $argT_{2}v(0)=\{0\}$ can be proved in a similar manner.\\

 To prove part (iii), we consider the optimal policy $a$ in node 1 service resource allocation. Since the service resources in node 1 is from the compact set $[0,a_{max}]$, the optimal policy $a$ in node 1 can be $0$, or $a_{max}$, or satisfies the following equation:
\begin{eqnarray*}\label{eq:nc1}
\frac{c'_{1}(a)}{\mu'_{1}(a)}=v(x)-v(x-e_{1}+e_{2}).
\end{eqnarray*}
We use the contradiction method. Assume that $a\in argT_{1}v(x)$ such that $a\in (0,a_{max})$ for all $x\in E$. For any $\varepsilon>0$, we have:
\begin{eqnarray*}
\lefteqn{T^{a+\varepsilon}_{1}v(x)-T^{a}_{1}v(x)}\\
&&=[\mu_{1}(a+\varepsilon)-\mu_{1}(a)][v(x-e_{1}+e_{2})-v(x)]+c_{1}(a+\varepsilon)-c_{1}(a)\geq 0,
\end{eqnarray*}
which implies that
\begin{eqnarray*}\label{eq:nc1}
 v(x)-v(x-e_{1}+e_{2})\leq\frac{c_{1}(a+\varepsilon)-c_{1}(a)}{\mu_{1}(a+\varepsilon)-\mu_{1}(a)}.
\end{eqnarray*}

Since the function $\frac{c_{1}(a)}{\mu_{1}(a)}$ is non-increasing, we get $\frac{c_{1}(a+\varepsilon)-c_{1}(a)}{\mu_{1}(a+\varepsilon)-\mu_{1}(a)}\leq\frac{c_{1}(a)}{\mu_{1}(a)} $, $v(x)-v(x-e_{1}+e_{2})\leq\frac{c_{1}(a)}{\mu_{1}(a)}$ which is a contradiction with the condition $\frac{c'_{1}(a)}{\mu'_{1}(a)}>\frac{c_{1}(a)}{\mu_{1}(a)}$. So there is no $a$ satisfying the above equation. That is, the optimal policy in node 1 is $ argT_{1}v(x)=\{0,a_{max}\}$. Thus, the optimal policy is a bang-bang control policy. The part (iii) for node 2 can be proved in a similar manner.

\end{pf}

\section{Conclusion}\label{sec5}
In this paper we have analysed the optimal server resources control of a tandem queueing system with two nodes. The controller can make a dynamic decision to allocate the service resource to each node at any decision epoch. Applying the dynamic programming to the model, we not only give some traditional properties of the relative value function and optimal policy, but also derive the condition under which the optimal policy is unique and bang-bang control occurs. In particular, we have provided the relationship between the two nodes' optimal policy, which can give the controller more information to manage the system.

From the above results there arise some interesting extensions of the model which we may study in the near future.

(i) One possible change is to consider a model where each node's service resource decision is dependent on the number of the customers in two queues. When the system state is $x=(x_{1},x_{2})$, the controller makes an action $f_{1}(x_{1},x_{2})=a\in A,f_{2}(x_{1},x_{2})=b\in B$. Although the analysis is difficult, we may get some another properties of the queue optimal policy. In our model the two nodes have their action sets. We can also study the further model in which the two nodes share the common server resources.

(ii) Another way to generalize the model is to consider some strategies in our model, such as the retrial, feedback and priority customers. The model may become more complex. Some other methods should be considered. In our model the customers arrive at the system according to a Poisson process and the service time of a customer is exponentially distributed. We can apply the embedded Markov chain and semi-Markov decision processes to consider the queueing system in which the service time of a customer is a general distribution.

(iii) In addition, the tandem queueing system with $n$ nodes is also worthy thinking about. Based on our model, we can study the optimal policy relationship between the two nodes.

\section*{Appendix A}
Property 4.1 (non-decreasingness)
\begin{pf}
To prove Property 4.1 (i), the proof is done by induction on $n$ in $v_{n}$. Define $v_{0}(x)=0$ for all state $x\in E$. This function obviously satisfies (i). Now, we assume that (i) holds for the function $v_{n}(x)$,$x\in E$ and some $n\in N$. We should prove that $v_{n+1}(x)$ satisfies the non-decreasing property as well. Then for $i=1$, we can get
\begin{eqnarray*}
\lefteqn{v_{n+1}(x+e_{1})-v_{n+1}(x)}\\
&&=\lambda[v_{n}(x+2e_{1})-v_{n}(x+e_{1})]+h_{1}+\sum_{i=1,2}T_{i}v_{n}(x+e_{1})-\sum_{i=1,2}T_{i}v_{n}(x).
\end{eqnarray*}
The second term of the right-hand side is obviously positive.

Let $(a\in argT_{1}v(x),b\in argT_{2}v_{n}(x))$ be an arbitrary optimal policy for node 1 and 2 in state $x$, respectively. Then
\begin{eqnarray*}
\lefteqn{\sum_{i=1,2}T_{i}v_{n}(x+e_{1})-\sum_{i=1,2}T_{i}v_{n}(x)}\\
&&\geq\mu_{1}(a)[v_{n}(x+e_{2})-v_{n}(x+e_{2}-e_{1})]\\
&&\  \  +\mu_{2}(b)[v_{n}(x-e_{2}+e_{1})-v_{n}(x-e_{2})]\\
&&\  \  +[\mu_{1}(a_{max})-\mu_{1}(a)+\mu_{2}(b_{max})-\mu_{2}(b)][v_{n}(x+e_{1})-v_{n}(x)]\\
&&\geq0,
\end{eqnarray*}
Therefore, Property 4.1 (i) holds by induction for any $n$, $v(x)$ is a nondecreasing function. Property 4.1 (i) for $i=2$ can be proved in a similar manner.

To prove Property 4.1 (ii), the proof is similar to the above one. Define $v_{0}(x)=0$ for all state $x\in E$. This function obviously satisfies the (ii). Now, we assume that (ii) holds for function $v_{n}(x)$, $x\in E$ and some $n\in N$. We should prove that $v_{n+1}(x)$ satisfies Property 4.1 (ii) as well.
\begin{eqnarray*}
\lefteqn{v_{n+1}(x-e_{1}+e_{2})-v_{n+1}(x-e_{2})}\\
&&=\lambda[v_{n}(x+e_{2})-v_{n}(x+e_{1}-e_{2})]+2h_{2}-h_{1}\\
&&\  \  +\sum_{i=1,2}T_{i}v_{n}(x-e_{1}+e_{2})-\sum_{i=1,2}T_{i}v_{n}(x-e_{2}).
\end{eqnarray*}
Since the condition $2h_{2}\geq h_{1}$ holds, the second term of the right-hand side is obviously positive.

Let $(a\in argT_{1}v(x-e_{2}),b\in argT_{2}v(x-e_{2}))$ be an arbitrary optimal policy for node 1 and 2 in state $x-e_{2}$, respectively. Then
\begin{eqnarray*}
\lefteqn{\sum_{i=1,2}T_{i}v_{n}(x-e_{1}+e_{2})-\sum_{i=1,2}T_{i}v_{n}(x-e_{2})}\\
&&\geq\mu_{1}(a)[v_{n}(x-2e_{1}+2e_{2})-v_{n}(x-e_{2})]\\
&&\  \  +\mu_{2}(b)[v_{n}(x-e_{1})-v_{n}(x-2e_{2})]\\
&&\  \  +[\mu_{1}(a_{max})-\mu_{1}(a)][v_{n}(x-e_{1}+e_{2})-v_{n}(x-e_{2})]\\
&&\  \  +[\mu_{2}(b_{max})-\mu_{2}(b)][v_{n}(x-e_{1}+e_{2})-v_{n}(x-e_{2})]\\
&&\geq0.
\end{eqnarray*}
Therefore, Property4.1 (ii) holds by induction for any $n$, we have $v(x-e_{1}+e_{2})\geq v(x-e_{2})$ for all $x=(x_{1},x_{2})\in E$ and $x_{1}\geq1,x_{2}\geq1$. Property 4.1 (iii) can be proved in a similar manner.
\end{pf}

\section*{Appendix B}
Property 4.2 (quasi-convexity) (i) and Theorem 1 (i)
\begin{pf}
 To prove Property 4.2 (i), we assume that Property 4.2 (i) for function $v_{n}(x)$, $x\in E$ and some $n\in N$ holds. Then we need to prove that Property 4.2 (i) for $n+1$ also holds. When $x=(x_{1},x_{2})\in E$ and $x_{2}\geq1$, we have
\begin{eqnarray*}
\lefteqn{v_{n+1}(x+e_{2})-2v_{n+1}(x)+v_{n+1}(x-e_{2})}\\
&&=\lambda[v_{n}(x+e_{2}+e_{1})-2v_{n}(x+e_{1})+v_{n}(x+e_{1}-e_{2})]\\
&&\  \  +\sum_{i=1,2}T_{i}v_{n}(x+e_{2})-2\sum_{i=1,2}T_{i}v_{n}(x)+\sum_{i=1,2}T_{i}v_{n}(x-e_{2})\\
&&\geq\sum_{i=1,2}T_{i}v_{n}(x+e_{2})-2\sum_{i=1,2}T_{i}v_{n}(x)+\sum_{i=1,2}T_{i}v_{n}(x-e_{2}).
\end{eqnarray*}

The inequality holds by the induction hypothesis. The optimal policy of node 1 is only dependent on the number of customers in node 1 and the state $x+e_{2}$, $x$, $x-e_{2}$ have the same first entry $x_{1}$. Hence, they have the same optimal policy in node 1. We assume that $a\in argT_{1}v(x+e_{2}),b_{1}\in argT_{2}v(x+e_{2})$, $a\in argT_{1}v(x-e_{2}),b_{2}\in argT_{2}v(x-e_{2})$. Therefore, we get
\begin{eqnarray*}
\lefteqn{\sum_{i=1,2}T_{i}v_{n}(x+e_{2})-2\sum_{i=1,2}T_{i}v_{n}(x)+\sum_{i=1,2}T_{i}v_{n}(x-e_{2})}\\
&&\geq\mu_{1}(a)[v_{n}(x-e_{1}+2e_{2})-2v_{n}(x-e_{1}+e_{2})+v_{n}(x-e_{1})]\\
&&\  \  +[\mu_{1}(a_{max})-\mu_{1}(a)][v_{n}(x+e_{2})-2v_{n}(x)+v_{n}(x-e_{2})]\\
&&\  \  +[\mu_{2}(b_{1})-\mu_{2}(b_{2})][v_{n}(x)-v_{n}(x-e_{2})]\\
&&\  \  +\mu_{2}(b_{2})[v_{n}(x)-2v_{n}(x-e_{2})+v_{n}(x-2e_{2})]\\
&&\  \  +[\mu_{2}(b_{max})-\mu_{2}(b_{1})][v_{n}(x+e_{2})-v_{n}(x)]\\
&&\  \  +[\mu_{2}(b_{max})-\mu_{2}(b_{2})][v_{n}(x-e_{2})-v_{n}(x)]\\
&&=\mu_{1}(a)[v_{n}(x-e_{1}+2e_{2})-2v_{n}(x-e_{1}+e_{2})+v_{n}(x-e_{1})]\\
&&\  \  +[\mu_{1}(a_{max})-\mu_{1}(a)][v_{n}(x+e_{2})-2v_{n}(x)+v_{n}(x-e_{2})]\\
&&\  \  +\mu_{2}(b_{2})[v_{n}(x)-2v_{n}(x-e_{2})+v_{n}(x-2e_{2})]\\
&&\  \  +[\mu_{2}(b_{max})-\mu_{2}(b_{1})][v_{n}(x+e_{2})-2v_{n}(x)+v_{n}(x-e_{2})]\\
&&\geq0.
\end{eqnarray*}

The first inequality follows by taking a potentially suboptimal action in the second term of $\sum_{i=1,2}T_{i}v_{n}(x+e_{2})-2\sum_{i=1,2}T_{i}v_{n}(x)+\sum_{i=1,2}T_{i}v_{n}(x-e_{2})$. The equality follows by rearranging the terms. The last inequality follows by the induction hypothesis. Hence, we have $ v(x+e_{2})-2v(x)+v(x-e_{2})\geq0$.

For Theorem 1 (i), let $(b_{1}\in argT_{2}v(x+e_{2}),b_{2}\in argT_{2}v(x))$ be an optimal policy for node 2 in states $x+e_{2}$, $x$, respectively. The proof is done by contradiction. Suppose that $b_{1}<b_{2}$, then
\begin{eqnarray*}
\lefteqn{T^{b_{1}}_{2}v(x)-T^{b_{2}}_{2}v(x)}\\
&&=[\mu_{2}(b_{2})-\mu_{2}(b_{1})][v(x)-v(x-e_{2})]-[c_{2}(b_{2})-c_{2}(b_{1})]\geq0.
\end{eqnarray*}
Since Property 4.1 (i) above and $\mu_{2}(b_{2})-\mu_{2}(b_{1})>0$ holds, we have
\begin{eqnarray*}
\lefteqn{T^{b_{1}}_{2}v(x+e_{2})-T^{b_{2}}_{2}v(x+e_{2})}\\
&&=[\mu_{2}(b_{2})-\mu_{2}(b_{1})][v(x+e_{2})-v(x)]-[c_{2}(b_{2})-c_{2}(b_{1})]\\
&&>[\mu_{2}(b_{2})-\mu_{2}(b_{1})][v(x)-v(x-e_{2})]-[c_{2}(b_{2})-c_{2}(b_{1})]\\
&&\geq0.
\end{eqnarray*}
However, this implies that $b_{1}$ is not an optimal policy for node 2 in state $x+e_{2}$. Hence $b_{1}\geq b_{2}$ .
\end{pf}

Property 4.2(quasi-convexity) (ii) and Theorem 1 (ii)

To prove Property 4.2 (ii), we assume that Property 4.2 (ii) holds for function $v_{n}(x)$, $x\in E$ and some $n\in N$. Then we need to prove that Property 4.2 (ii) for $n+1$ also holds. When $x=(x_{1},x_{2})\in E$ and $x_{1}\geq1,x_{2}\geq1$, we have
\begin{eqnarray*}
\lefteqn{v_{n+1}(x+e_{1}-e_{2})-2v_{n+1}(x)+v_{n+1}(x-e_{1}+e_{2})}\\
&&=\lambda[v_{n}(x+2e_{1}-e_{2})-2v_{n}(x+e_{1})+v_{n}(x+e_{2})]\\
&&\  \  +\sum_{i=1,2}T_{i}v_{n}(x+e_{1}-e_{2})-2\sum_{i=1,2}T_{i}v_{n}(x)+\sum_{i=1,2}T_{i}v_{n}(x-e_{1}+e_{2})\\
&&\geq\sum_{i=1,2}T_{i}v_{n}(x+e_{1}-e_{2})-2\sum_{i=1,2}T_{i}v_{n}(x)+\sum_{i=1,2}T_{i}v_{n}(x-e_{1}+e_{2})\\
&&=T_{1}v_{n}(x+e_{1}-e_{2})-2T_{1}v_{n}(x)+T_{1}v_{n}(x-e_{1}+e_{2})\\
&&\  \  +T_{2}v_{n}(x+e_{1}-e_{2})-2T_{2}v_{n}(x)+T_{2}v_{n}(x-e_{1}+e_{2}).
\end{eqnarray*}

The inequality above holds by the induction hypothesis. Now, we assume that $a_{1}\in argT_{1}v(x+e_{1}-e_{2}),b_{1}\in argT_{2}v(x+e_{1}-e_{2})$, $a_{2}\in argT_{1}v(x-e_{1}+e_{2}),b_{2}\in argT_{2}v(x-e_{1}+e_{2})$. Then, we get
\begin{eqnarray*}
\lefteqn{T_{1}v_{n}(x+e_{1}-e_{2})-2T_{1}v_{n}(x)+T_{1}v_{n}(x-e_{1}+e_{2})}\\
&&\geq\mu_{1}(a_{1})[v_{n}(x)-v_{n}(x-e_{1}+e_{2})]\\
&&\  \  +\mu_{1}(a_{2})[v_{n}(x-2e_{1}+2e_{2})-v_{n}(x-e_{1}+e_{2})]\\
&&\  \  +[\mu_{1}(a_{max})-\mu_{1}(a_{1})][v_{n}(x+e_{1}-e_{2})-v_{n}(x)]\\
&&\  \  +[\mu_{1}(a_{max})-\mu_{1}(a_{2})][v_{n}(x-e_{1}+e_{2})-v_{n}(x)]\\
&&=\mu_{1}(a_{2})[v_{n}(x-2e_{1}+2e_{2})-2v_{n}(x-e_{1}+e_{2})+v_{n}(x)]\\
&&\  \  +[\mu_{1}(a_{max})-\mu_{1}(a_{1})][v_{n}(x+e_{1}-e_{2})-2v_{n}(x)+v_{n}(x-e_{1}+e_{2})]\\
&&\geq0.
\end{eqnarray*}

The first inequality follows by taking a potentially suboptimal action in the second term of the operator $T_{1}v_{n}(x+e_{1}-e_{2})-2T_{1}v_{n}(x)+T_{1}v_{n}(x-e_{1}+e_{2})$. The equality follows by rearranging the terms. The last inequality follows by the induction hypothesis.
\begin{eqnarray*}
\lefteqn{T_{2}v_{n}(x+e_{1}-e_{2})-2T_{2}v_{n}(x)+T_{2}v_{n}(x-e_{1}+e_{2})}\\
&&\geq\mu_{2}(b_{1})[v_{n}(x+e_{1}-2e_{2})-v_{n}(x-e_{2})]\\
&&\  \  +\mu_{2}(b_{2})[v_{n}(x-e_{1})-v_{n}(x-e_{2})]\\
&&\  \  +[\mu_{2}(b_{max})-\mu_{2}(b_{1})][v_{n}(x+e_{1}+e_{2})-v_{n}(x)]\\
&&\  \  +[\mu_{2}(b_{max})-\mu_{2}(b_{2})][v_{n}(x-e_{1}+e_{2})-v_{n}(x)]\\
&&=\mu_{2}(b_{2})[v_{n}(x+e_{1}-2e_{2})-2v_{n}(x-e_{2})+v_{n}(x-e_{1})]\\
&&\  \  +[\mu_{2}(b_{max})-\mu_{2}(b_{2})][v_{n}(x+e_{1}-e_{2})-2v_{n}(x)+v_{n}(x-e_{1}+e_{2})]\\
&&\  \  +[\mu_{2}(b_{1})-\mu_{2}(b_{2})][v_{n}(x+e_{1}-2e_{2})-v_{n}(x+e_{1}-e_{2})]\\
&&\geq0.
\end{eqnarray*}

The first inequality follows by taking a potentially suboptimal action in the second term of the operator above. The equality follows by rearranging the terms. The last one follows by the induction hypothesis and because of Theorem 1 (i), we know that $b_{1}\leq b_{2}$. So that we have $\mu_{2}(b_{1})-\mu_{2}(b_{2})\leq0$. From the Property 4.1, we know that $v_{n}(x+e_{1}-2e_{2})-v_{n}(x+e_{1}-e_{2})\leq0$. Thus, we derive that $[\mu_{2}(b_{1})-\mu_{2}(b_{2})][v_{n}(x+e_{1}-2e_{2})-v_{n}(x+e_{1}-e_{2})]\geq0$. Therefore, the last inequality is taken.

For Theorem 1 (ii), let $(a_{1}\in argT_{1}v(x+e_{1}-e_{2}),a_{2}\in argT_{1}v(x))$ be an optimal policy for node 2 in states $x+e_{1}-e_{2}$, $x$, respectively. The proof is done by contradiction. Suppose that $a_{1}<a_{2}$, then
\begin{eqnarray*}
\lefteqn{T^{a_{1}}_{1}v(x)-T^{a_{2}}_{1}v(x)}\\
&&=[\mu_{1}(a_{2})-\mu_{1}(a_{1})][v(x-e_{1}+e_{2})-v(x)]-[c_{1}(a_{2})-c_{1}(a_{1})]\\
&&\geq0.
\end{eqnarray*}
From Property 4.1 (ii) above and $\mu_{1}(a_{2})-\mu_{1}(a_{1})>0$, we have
\begin{eqnarray*}
\lefteqn{T^{a_{1}}_{1}v(x+e_{1}-e_{2})-T^{a_{2}}_{1}v(x+e_{1}-e_{2})}\\
&&=[\mu_{1}(a_{2})-\mu_{1}(a_{1})][v(x)-v(x+e_{1}-e_{2})]-[c_{1}(a_{2})-c_{1}(a_{1})]\\
&&\geq[\mu_{1}(a_{2})-\mu_{1}(a_{1})][v(x-e_{1}+e_{2})-v(x)]-[c_{1}(a_{2})-c_{1}(a_{1})]\\
&&\geq0.
\end{eqnarray*}
However, this implies that $a_{1}$ is not an optimal policy for node 1 in state $x+e_{1}-e_{2}$. Hence $a_{1}\geq a_{2}$.

Since the optimal policy of node 1 is dependent only on the number of customers in node 1, and the states $x+e_{1}$, $x+e_{1}-e_{2}$ have the same first entry $x_{1}+1$. So they have the same optimal policy $a_{1}$ in node 1, i.e., $a_{1}\in argT_{1}v(x+e_{1})$. Thus we get that if $ a_{1}\in argT_{1}v(x+e_{1}),a_{2}\in argT_{1}v(x)$ hold, then we have $a_{1}\geq a_{2}$ for all $x=(x_{1},x_{2})\in E$.

\end{document}